\documentclass[reqno]{amsart}
\usepackage{amsmath,amsthm,amscd,amssymb,amsfonts, amsbsy}
\usepackage{latexsym, color, enumerate}
\usepackage{pxfonts}

\numberwithin{equation}{section}

\theoremstyle{plain}
\newtheorem{theorem}[equation]{Theorem}

\theoremstyle{definition}

\theoremstyle{remark}
\newtheorem{remark}[equation]{Remark}

\newcommand{\dv}{\operatorname{div}}
\newcommand{\R}{\operatorname{Re}}

\newcommand{\diam}{\operatorname{diam}}

\newcommand{\mysection}[1]{\section{#1}
\setcounter{equation}{0}}

\newcommand{\bC}{\mathbb{C}}

\newcommand{\bR}{\mathbb{R}}

\providecommand{\set}[1]{\{#1\}}

\providecommand{\abs}[1]{\lvert#1\rvert}
\providecommand{\Abs}[1]{\left\lvert#1\right\rvert}

\providecommand{\norm}[1]{\lVert#1\rVert}

\renewcommand{\qedsymbol}{$\blacksquare$}

\begin{document}
\title[Elliptic equation arising from photo-acoustic imaging]
{On an elliptic equation arising from photo-acoustic imaging in inhomogeneous media}

\author[H. Ammari]{Habib Ammari}
\address[H. Ammari]{Department of Mathematics and Applications, Ecole Normale Sup\'erieure,
45 Rue d'Ulm, 75005 Paris, France}
\email{habib.ammari@ens.fr}
\thanks{H. Ammari was partially supported by the
ERC Advanced
Grant Project MULTIMOD--267184.}

\author[H. Dong]{Hongjie Dong}
\address[H. Dong]{Division of Applied Mathematics, Brown University,
182 George Street, Providence, RI 02912, United States of America}
\email{hdong@brown.edu}
\thanks{H. Dong was partially supported by the NSF under agreement DMS-1056737.}

\author[H. Kang]{Hyeonbae Kang}
\address[H. Kang]{Department of Mathematics, Inha University, Incheon 402-751, Republic of Korea}
\email{hbkang@inha.ac.kr}
\thanks{H. Kang was partially supported by the Korean Ministry of Education, Sciences and Technology through NRF grant No. 2010-0017532.}

\author[S. Kim]{Seick Kim}
\address[S. Kim]{Department of Mathematics, Yonsei University, 50 Yonsei-ro, Seodaemun-gu, Seoul 120-749, Republic of Korea}
\email{kimseick@yonsei.ac.kr}
\thanks{S. Kim is partially supported by NRF-2012R1A1A2040411.}

\subjclass[2010]{Primary 35B45, 35B65; Secondary 35J08}

\keywords{De Giorgi estimate, Green function, measurable coefficients.}

\begin{abstract}
We study an elliptic equation with measurable coefficients arising from photo-acoustic imaging in inhomogeneous media.
We establish H\"older continuity of weak solutions and obtain pointwise bounds for Green's functions subject to Dirichlet or Neumann condition.
\end{abstract}

\maketitle

\mysection{Introduction and main results}
We consider an operator $L$ of the form
\[
L u = -\dv (\gamma (x) \nabla u) + ik u
\]
acting on a complex-valued function $u$ defined on a domain $\Omega \subset \bR^n$ ($n\ge 3$).
Here,  $i=\sqrt{-1}$,  $k$ is a nonzero constant, and $\gamma(x)$ is a real $n\times n$ matrix valued function satisfying
\begin{equation}		\label{eq012ho}
\nu \abs{\xi}^2 \le  \R \bar{\xi}^\top \gamma(x) \xi,  \quad  \forall \xi \in \bC^n,\quad \forall x \in \Omega;\qquad
\norm{\gamma}_{L^\infty(\Omega)} \le \nu^{-1},
\end{equation}
for some $\nu\in (0,1]$.
The adjoint operator $L^*$ is given by
\[
L^* u = -\dv (\gamma(x)^\top \nabla u) - i k u.
\]
Note that the condition \eqref{eq012ho} is invariant if we replace $\gamma$ by $\gamma^\top$.
Therefore, even if $L$ is not necessarily a symmetric operator, $L$ and $L^*$ are operators of the same type.

The operator $L$ arises from quantitative photo-acoustic imaging in inhomogeneous media; see \cite{ABJK, ABJW}. The purpose of quantitative photo-acoustic imaging is to reconstruct the optical absorption coefficient of an optical anomaly from the absorbed energy density. The absorbed energy density is obtained from boundary measurements of the pressure induced by the thermo-elastic effect; see \cite{wang}. The reconstruction problem is challenging. One promising method is implemented in \cite{ABJK2} in the framework of diametrically small optical anomalies. It uses asymptotic expansions of the energy density in terms of the characteristic size of the optical anomalies. Such expansions are based on estimates of the Neumann function of the operator $L$ under the assumption that $\gamma$ is uniformly continuous in $\Omega$. This note aims at removing such regularity assumption and extending the results of \cite{AKK12} to real biological tissues.

We say that the function $u \in H^1(\Omega)$ is a weak solution of $L u = f$,  $f\in L^{2n/(n+2)}(\Omega)$, if it satisfies for any $\varphi \in H^1_0(\Omega)$ the identity
\begin{equation}		\label{eq013jp}
\int_\Omega  \left(\nabla \bar{\varphi}^\top \gamma \nabla u  + i k u \bar{\varphi} \right) =\int_\Omega f \bar{\varphi}.
\end{equation}
In this short article, we prove that any weak solution $u$ of $Lu=0$ satisfies the following De Giorgi-Moser-Nash type estimate in $B_r(x_0) \subset \Omega$:
\[
[u]_{\alpha; B_{r/2}(x_0)} \le C  r^{-n/2 {\color{blue}-} \alpha} \norm{u}_{L^2(B_r(x_0))},
\]
where $[u]_{\alpha; B_r(x_0)}$ is the standard H\"older semi-norm of $u$ over an open ball of radius $r$ centered at $x_0$.
We emphasize that the constant $C$ in the above estimate is independent of $k$ and depends only on $n$ and $\nu$.
Since $u$ is complex valued and its real and imaginary parts are intertwined in the equation through the parameter $k$, this result is not immediate from the classical De Giorgi-Moser-Nash theory.
We also establish similar estimates near the boundary for weak solutions satisfying zero Dirichlet or Neumann boundary condition.
Then by utilizing these estimates, we construct Dirichlet and Neumann functions and derive the pointwise estimates for them; we show that they are bounded in modulus by $C \abs{x-y}^{2-n}$.
By Dirichlet's function we mean a function $G(x,y)$ satisfying
\[
\left\{
\begin{aligned}
L G(\cdot,y) &=\delta_y \quad \text{in }\;\Omega,\\
G(\cdot,y)&= 0 \quad \text{on }\;\partial\Omega,
\end{aligned}
\right.
\]
and similarly, by Neumann function, we mean a function $N(x,y)$ satisfying
\[
\left\{
\begin{aligned}
L N(\cdot,y) &=\delta_y \quad \text{in }\;\Omega,\\
\gamma \nabla N(\cdot,y)\cdot n &= 0 \quad \text{on }\;\partial\Omega.
\end{aligned}
\right.
\]
The above definitions for Dirichlet and Neumann functions should be interpreted in a weak sense; see \cite[Sec.~2.2]{AKK12} for details.

Now, let us state our main results more precisely.
\begin{theorem}			\label{thm1}
Let $u \in H^1(\Omega)$ be a weak solution of $Lu=f$.
\begin{enumerate}[(i)]
\item
Let $B_r=B_r(x_0) \subset \Omega$ and $f \in L^p(B_r)$ for some $p>n/2$.
Then, we have
\begin{equation}		\label{eq014gh}
[u]_{\alpha; B_{r/2}} \le C r^{-\alpha} \left( r^{-n/2} \norm{u}_{L^2(B_r)} + r^{2-n/p} \norm{f}_{L^p(B_r)}\right),
\end{equation}
where $\alpha\in (0,1)$ and $C>0$ are constants depending only on $n$, $\nu$, and $p$.
\item
Suppose $u=0$ on $\partial\Omega$ in the sense that $u \in H^1_0(\Omega)$.
Let $\Omega_r=\Omega_r(x_0)=\Omega \cap B_r(x_0)$ with $x_0\in \bar \Omega$ and $f \in L^p(\Omega_r)$ for some $p>n/2$.
Then, we have
\begin{equation}		\label{eq014hb}
\norm{u}_{L^\infty(\Omega_{r/2})} \le C \left( r^{-n/2} \norm{u}_{L^2(\Omega_r)} + r^{2-n/p} \norm{f}_{L^p(\Omega_r)}\right)
\end{equation}
for any $0<r<\diam\Omega$, where $C$ is a constant depending only on $n$, $\nu$, and $p$.

\item
Suppose that $\Omega$ is a bounded Lipschitz domain and  $\partial u /\partial \nu =0$ on $\partial \Omega$ in the sense that for any $\varphi \in H^1(\Omega)$ the identity \eqref{eq013jp} holds.
Then, there is a constant $C$ depending only on $n$, $\nu$, $p$, and $\Omega$ such that the estimate \eqref{eq014hb} is valid for any $0<r<\diam\Omega$.
\end{enumerate}
\end{theorem}

\begin{remark}
It should be clear from the proof that the same conclusions of the preceding theorem hold true for any weak solution $u \in H^1(\Omega)$ of $L^* u = f$.
\end{remark}

\begin{remark}
In view of the proof of Theorem \ref{thm1} below, by using the boundary De Giorgi-Moser-Nash theorem, it is easily seen that under the assumption of Theorem~\ref{thm1} ii) and a uniform exterior cone condition on $\Omega$,
we have the following H\"older estimate:
\[
[u]_{\alpha; \Omega_{r/2}} \le C r^{-\alpha} \left( r^{-n/2} \norm{u}_{L^2(\Omega_r)} + r^{2-n/p} \norm{f}_{L^p(\Omega_r)}\right),
\]
where $\alpha\in (0,1)$ and $C>0$ are constants depending only on $n$, $\nu$, $p$ and the parameters in a uniform exterior cone condition.
The same estimate also holds under the assumption of Theorem~\ref{thm1} iii).
\end{remark}

\begin{theorem}		\label{thm2}
(Dirichlet's function)
There exists a unique Dirichlet function $G(x,y)$ such that
\begin{equation}		\label{eq015ba}
\abs{G(x,y)} \le C \abs{x-y}^{2-n},\quad \forall x,y \in \Omega,\quad x\neq y,
\end{equation}
where $C$ is a constant depending only on $n$ and $\nu$.
In particular, when $\Omega=\bR^n$, we have the fundamental solution with the pointwise estimate \eqref{eq015ba}.
\end{theorem}

\begin{theorem}		\label{thm3}
(Neumann function)
Let $\Omega$ be a bounded Lipschitz domain.
There exists a unique Neumann function $N(x,y)$ such that
\begin{equation}		\label{eq017hg}
\abs{N(x,y)} \le C  \max(1,\abs{k}^{-1/2})\; \abs{x-y}^{2-n},\quad \forall x,y \in \Omega,\quad x\neq y,
\end{equation}
where $C$ is a constant depending only on $n$, $\nu$, and $\Omega$.
\end{theorem}

We note that a similar result regarding Neumann function $N(x,y)$ was established in \cite{AKK12} under the stronger assumption that $\gamma$ is uniformly continuous in $\Omega$ and $\Omega$ is a $C^1$-domain.
Therefore, our result is quite an improvement of that in \cite{AKK12}.
We emphasize, however, that in \cite{AKK12} estimates of derivatives of the Neumann functions are also derived under additional smoothness assumptions on the coefficients and the domain.

Finally, one may wish to consider the case when $k$ is not a constant but a variable coefficients $k(x)$ and  include lower order terms in $L$.
Even in this case, our main estimates are independent of  the size of $k(x)$ but on the size of its ``oscillation''.
We shall illustrate how to generalize our results to this case in Section \ref{sec2.4}.

\mysection{Proof of main theorems} \label{proofs}
\subsection{Proof of Theorem~\ref{thm1}}
The idea of the proof is to add a new variable and rewrite the equation as a parabolic equation without zero-order term, so that $k$ does not play a role in our estimates.
Let $v(x,t)=e^{ikt}u(x)$. Then $v$ becomes a weak solution of the parabolic equation
\[
v_t-\dv(\gamma(x) \nabla v)=e^{ikt}f(x)\quad\text{in }\;\Omega\times (-\infty,\infty).
\]
Note that the real and imaginary parts of $v$ are solutions of a real scalar equation
\[
w_t-\dv(\gamma(x) \nabla w)=g(x,t)\quad\text{in }\;\Omega\times (-\infty,\infty)
\]
with $g(x,t)=\cos(kt) f(x)$ or $g(x,t)=\sin(kt) f(x)$.
Denote $Q_r=B_r \times (-r^2,0)$.
Clearly, we have
\[
\norm{g}_{L^{p, \infty}_{x,t}(Q_r)} \le \norm{f}_{L^p(B_r)}
\]
and thus by the De Giorgi-Moser-Nash theorem applied to the real and imaginary parts of $v$, we have the parabolic H\"older estimate
\[
[v]_{C^{\alpha,\alpha/2}(Q_{r/2})} \le C r^{-\alpha} \left(r^{-(n+2)/2} \norm{v}_{L^2(Q_r)} + r^{2-n/p} \norm{g}_{L^\infty_tL^p_x(B_r)}\right),
\]
where $C$ depends only on $n$, $\nu$, and $p$.

Then for any $x, y \in B_{r/2}(x_0)$ with $x\neq y$, we have
\begin{align*}
\frac{\abs{u(x)-u(y)}}{\abs{x-y}^\alpha}&=\frac{\abs{v(x,t)-v(y,t)}}{\abs{x-y}^\alpha}\\
&\le C r^{-\alpha} \left( r^{-(n+2)/2} \norm{v}_{L^2(Q_r)}+ r^{2-n/p} \norm{g}_{L^{p,\infty}_{x,t}(Q_r)}\right) \\
&=C r^{-\alpha} \left( r^{-n/2} \norm{u}_{L^2(B_r)}+ r^{2-n/p} \norm{f}_{L^p(B_r)}\right).
\end{align*}
Therefore, we have the estimate \eqref{eq014gh}.
This completes the proof for part (i).

Next suppose that $u \in H^1_0(\Omega)$.
Then $v(x,t)=e^{ikt}u(x)$ vanishes on $\partial\Omega\times (-\infty,\infty)$.
Denote $U_r=\Omega_r \times (-r^2, 0)$.
By the local maximum principle of De Giorgi and Moser (see, e.g., \cite[Ch. III, \S~7]{LSU}) applied to both real and imaginary parts of $v$, we get
\begin{equation}		\label{eq200dgm}
\norm{v}_{L^{\infty}(U_{r/2})} \le C  \left(r^{-(n+2)/2} \norm{v}_{L^2(U_r)} + r^{2-n/p} \norm{g}_{L^{p,\infty}_{x,t}(U_r)}\right),
\end{equation}
where $C$ depends only on $n$, $\nu$, and $p$.
Then as in part (i), we get \eqref{eq014hb} from the above inequality.
This completes the proof for part (ii).
The proof for part (iii) is similar.
The assumption that $\Omega$ is a Lipschitz domain makes the (multiplicative) Sobolev inequality available to us and we can still apply De Girogi or Moser's method as mentioned above to get the local boundedness estimate \eqref{eq200dgm}.
The constant $C$ at this time, however, additionally depends on the constant of the Sobolev inequality and diameter of $\Omega$, and thus depends on $\Omega$; see \cite[\S~5.3.1]{CK14}.
\hfill\qedsymbol

\subsection{Proof of Theorem~\ref{thm2}}
First, we establish unique solvability of the problem
\begin{equation}		\label{eq202sol}
L u = f \quad\text{in }\;\Omega,\quad f \in L^{2n/(n+2)}(\Omega)
\end{equation}
in the space $H^1_0(\Omega)$.
For any $R>0$, denote $\Omega_R=\Omega \cap B_R(0)$ and consider the Hilbert space $H^1_0(\Omega_R)$ equipped with inner product
\[
(u,v)=\int_{\Omega_R} \nabla u \cdot \nabla \bar{v}.
\]
We define the sesquilinear form $B: H^1_0(\Omega_R) \times H^1_0(\Omega_R) \to \bC$ as
\begin{equation*}                  
B[u,v]:=\int_{\Omega_R} (\nabla \bar{v}^\top \gamma \nabla u   + ik u \bar v).
\end{equation*}
It is easy to check that $B$ satisfies the following:
\begin{enumerate}[i)]
\item
there exists $C=C(n,\nu, k, R)$ such that
\[
\abs{B[u,v]} \le C \norm{u}\, \norm{v}  ,\quad \forall u, v \in H^1_0(\Omega_R);
\]
\item
\[
\nu\norm{u}^2 \le \R B[u,u],\quad \forall u \in H^1_0(\Omega_R).
\]
\end{enumerate}
Let $F: H^1_0(\Omega_R) \to \bC$ be a bounded anti-linear functional; i.e.,
\[
F(v+\lambda w)= F(v)+ \bar{\lambda} F(w),\qquad \forall v, w \in H^1_0(\Omega_R), \quad \forall \lambda \in \bC
\]
and
\[
\norm{F}_*:=\sup_{\norm{v} \le 1} \abs{F(v)} <\infty.
\]
Then by the Lax-Milgram lemma (in the complex case), there exists a unique $u \in H^1_0(\Omega_R)$ such that
$\norm{u} \le \nu^{-1} \norm{F}_*$ and
\[
B[u,v]= F(v), \quad \forall v\in H^1_0(\Omega_R).
\]
In particular, for any $f\in L^{2n/(n+2)}(\Omega)$, there is a unique $u^{(R)} \in H^1_0(\Omega_R)$ such that
\begin{equation}		\label{eq2.05no}
\int_{\Omega_R} (\nabla \bar{v}^\top \gamma \nabla u^{(R)}   + ik u^{(R)} \bar v) = \int_{\Omega_R} f\bar{v},\quad \forall v\in H^1_0(\Omega_R)
\end{equation}
and $u^{(R)}$ satisfies the estimate
\begin{equation*}		
\norm{\nabla u^{(R)}}_{L^2(\Omega_R)} \le C(n) \nu^{-1} \norm{f}_{L^{2n/(n+2)}(\Omega)}.
\end{equation*}
By taking $v=u^{(R)}$ in \eqref{eq2.05no} and using Cauchy's inequality and the Sobolev embedding theorem, we get
\[
\sqrt{\abs{k}}\; \norm{u^{(R)}}_{L^2(\Omega_R)} \le C(n, \nu) \norm{f}_{L^{2n/(n+2)}(\Omega)}.
\]
Therefore, there exist a sequence $\set{u_j}$ and a function $u$ in $H^1_0(\Omega)$
such that $u_j \rightharpoonup u$ weakly in $H^1(\Omega)$, and $u$ satisfies
\begin{equation}
                            \label{eq11.57}
\norm{\nabla u}_{L^2(\Omega)} \le C(n) \nu^{-1} \norm{f}_{L^{2n/(n+2)}(\Omega)},\quad
\sqrt{\abs{k}}\; \norm{u}_{L^2(\Omega)} \le C(n, \nu) \norm{f}_{L^{2n/(n+2)}(\Omega)}.
\end{equation}
Thus by taking the limit in \eqref{eq2.05no}, we have
\begin{equation}
                                    \label{eq11.59}
\int_{\Omega} (\nabla \bar{v}^\top \gamma \nabla u   + ik u \bar v) = \int_{\Omega} f\bar{v},\quad \forall v\in C_0^\infty(\Omega).
\end{equation}
By \eqref{eq11.57} and a density argument, we get \eqref{eq11.59} for any $v\in H^1_0(\Omega)$.
This establishes existence of a solution in the space $H^1_0(\Omega)$.
In the case when $f=0$, the above identity implies that if $u\in H^1_0(\Omega)$ is a weak solution of $Lu=0$, then $u=0$.
This establishes uniqueness too.

Having unique solvability of the problem \eqref{eq202sol} and the part i) and ii) of Theorem~\ref{thm1}, one can construct Dirichlet's function and obtain the pointwise bound \eqref{eq015ba} by following the same method as in \cite{HK07} and \cite{KK10}.
The uniqueness of Dirichlet's function is a mere consequence of the fact that $u=0$ is a unique weak solution of $Lu=0$ in the function class of $H^1_0(\Omega)$.
\hfill \qedsymbol

\subsection{Proof of Theorem~\ref{thm3}}
The theorem follows from replicating the proof of \cite[Theorem~2.3]{AKK12} using part (i) and (iii) of Theorem~\ref{thm1}.
It is assumed that $k \ge k_0>0$ in \cite{AKK12}.
However, by replacing the constant $C=C(n, \nu, k_0)$ in the estimate \cite[(2.29)]{AKK12} by $\max(1, \abs{k}^{-1/2})\, C(n, \nu)$, we dispense with introducing the lower bound $k_0$.
\hfill \qedsymbol

\subsection{A generalization}        \label{sec2.4}

In this subsection, we show how to extend our results to operators with lower-order terms when $k$ is not a constant but a variable coefficients $k(x)$.
We consider an operator $L$ of the form
\begin{equation}		\label{eq:lot}
Lu= -\dv(\gamma(x) \nabla u+\tilde b(x) u) + b(x) \cdot \nabla u +  c(x) u + ik(x) u,
\end{equation}
where $b=(b^1,\ldots, b^n)^\top$, $\tilde{b}=(\tilde{b}^1,\ldots, \tilde{b}^n)^\top$, and $k(x)$ are real valued functions defined on $\Omega \subset \bR^n$.
We assume that $b$, $\tilde b$, and $c$ are bounded functions so that
\[
\norm{b}_\infty+\norm{\tilde{b}}_\infty+ \norm{c}_\infty \le M.
\]
When $k(x)\equiv k \in \bR$, by the same proof, the estimate \eqref{eq014gh} in Theorem~\ref{thm1} (i) remains valid for all $r \in (0, r_o]$, where $r_o>0$ is any fixed number, but the constant $C$ now depends on $\nu$, $n$, $p$, $r_o$, and $M$; but it is still independent of $k$.

Next we treat the case with general variable coefficient $k(x)$.
We assume that $k(x)$ satisfies
\[
\sup_{r \in (0, r_o]}\, r^{2} \left(\fint_{B_r} \Abs{k(x)-\bar k_{r}}^q\,dx\right)^{1/q}\le \kappa_o, \quad\text{where}\quad \bar k_{r} = \fint_{B_r} k(x)\,dx,
\]
for some $r_o>0$, $q>n/2$, and $\kappa_o \ge 0$.
This condition is satisfied, for instance, when $k(x)$ is uniformly continuous or is a function of bounded mean oscillations (BMO).
We rewrite the equation $Lu=f$ into
\[
-\dv(\gamma(x) \nabla u) + b(x) \cdot \nabla u +  c(x) u + i \bar k_{r} u=f- i(k(x)-\bar k_{r}) u.
\]
Choose $p_1\in (n/2,q)$ and $p_2\in (n/2,\infty)$ such that $1/p_1=1/q+1/p_2$.
By the previous statement regarding an operator $L$ with $k(x) \equiv k \in \bR$, we have
\begin{align*}
[u]_{\alpha; B_{r/2}} &\le C r^{-\alpha} \left( r^{-n/2} \norm{u}_{L^2(B_r)} + r^{2-n/p} \norm{f}_{L^p(B_r)}+r^{2-n/p_1}\norm{(k-\bar k_{r}) u}_{L^{p_1}(B_r)}\right)\\
&\le C r^{-\alpha} \left( r^{-n/2} \norm{u}_{L^2(B_r)} + r^{2-n/p} \norm{f}_{L^p(B_r)}+\kappa_o r^{-n/p_2}\norm{u}_{L^{p_2}(B_r)}\right),
\end{align*}
where we used H\"older's inequality in the second line.
Estimate \eqref{eq014gh} then follows from the inequality above and a standard iteration argument (cf. \cite[pp. 80--82]{Giaq}).
Of course, now the constant $C$ also depends on $q$ and $\kappa_o$ in addition to  $\nu$, $n$, $p$, $r_o$, and $M$, but is still independent of the size of $k(x)$.

Similar adjustments are needed for Theorem~\ref{thm1} (ii) and (iii).
Theorem~\ref{thm2} remains valid when $\Omega$ is bounded and we have coercivity for the sequilinear form associated with \eqref{eq:lot}; i.e., under these assumptions, there exists a unique Dirichlet function $G(x,y)$ such that
\[
\abs{G(x,y)} \le C \abs{x-y}^{2-n},
\]
where $C$ is a constant depending on $\nu$, $n$, $M$, $\kappa_o$, and $\diam \Omega$.
We have the required coercivity, for example, when $c(x) \ge c_o$ for some $c_o$ large enough.
Under an additional assumption that $k(x)$ is bounded away from zero, Theorem~\ref{thm3} also remains true if the estimate \eqref{eq017hg} is replaced by
\[
\abs{N(x,y)} \le C \max \left(1, (\inf \,\abs{k})^{-1/2}\right)\, \abs{x-y}^{2-n},
\]
where the constant $C$ depends additionally on $M$ and $\kappa_o$.

\begin{remark}
As in the classical De Giorgi--Nash--Moser theory, it is possible to relax assumptions on $b$, $\tilde{b}$, and $c$ so that they belong to certain Lebesgue or Morrey spaces other than the $L^\infty$ space.
We leave the details to the interested reader.
\end{remark}


\end{document}